\documentclass{article}
\usepackage{amssymb}
\usepackage{amsmath}
\usepackage{amsthm}

\newtheorem{theorem}{Theorem}
\newtheorem{lemma}{Lemma}[section]
\newtheorem{propo}{Proposition}
\newtheorem{question}{Question}[section]

\newtheorem{defin}{Definition}[section]
\newtheorem{remark}{Remark}
\newtheorem{problem}{Problem}

\newcommand{\bF}{\mathbb{F}}

\newcommand{\F}{\bF}

\def\proof{\smallskip\noindent{\it Proof.} }

\newcommand{\cM}{\mathcal{M}}

\def\vi{\vskip 0.1in \noindent}
\begin{document}

\title{The Space of Actions, Partition Metric and Combinatorial Rigidity}
\author{Mikl\'{o}s Ab\'{e}rt and G\'{a}bor Elek\thanks{%
AMS Subject Classification: 37A15 \thinspace\ The first author was partially supported by the ERC Consolidator Grant Invgrogra 648017, OTKA Grant 67867, NK 78439 and TAMOP 4.2.1/B-09/1/KMR-2010-003}}
\maketitle

\begin{abstract}
\noindent
We introduce a natural pseudometric on the space of actions of d-generated groups. In this pseudometric, the zero classes correspond to the weak equivalence classes defined by Kechris, and the metric identification is compact.
We achieve this by employing symbolic dynamics and an ultraproduct construction which also facilitates the extension of our results to unitary representations. 
As a byproduct, we show that the weak equivalence class of every
free non-amenable action contains an action that satisfies the measurable
von Neumann problem.

\end{abstract}

\section{Introduction}

Let $(X,\mathcal{B},\mu )$ be a non-atomic standard Borel probability space.
 Unless
otherwise specified, we assume $(X,\mathcal{B},\mu )$ is the unit
interval equipped with the usual Borel subsets and Lebesgue measure. An \emph{%
automorphism} of $(X,\mathcal{B},\mu )$ is defined as a measure preserving
Borel isomorphism of $(X,\mathcal{B},\mu )$. We identify two isomorphisms if
they act the same way up to a nullset. We define $\mathrm{Aut}(X,\mathcal{B},\mu )$
as the group of equivalence classes of automorphisms of $(X,\mathcal{B},\mu )$. Using a dense
sequence $\{A_{n}\}^\infty_{n=1}$ in the measure algebra $\cM(X,\mu)$ (see Section \ref{sect2}), one can
define a metric on $\mathrm{Aut}(X,\mathcal{B},\mu )$ as 
\begin{equation*}
\delta (\phi,\psi)=\sum\limits_{n}2^{-n}\mu \left( \phi(A_{n})\vartriangle
\psi(A_{n})\right)
\end{equation*}%
where $\vartriangle $ denotes symmetric difference. The metric $\delta $
induces the weak topology on $\mathrm{Aut}(X,\mathcal{B},\mu )$.
The metric space $(\mathrm{Aut}(X,\mathcal{B},\mu),\delta)$ is separable and
complete, making the group of automorphisms a Polish group. \vi
For a finite alphabet $S$, let 
$\mathrm{A}(S)=\mathrm{Aut}(X,\mathcal{B},\mu
)^{S}$ denote the set of maps from $S$ to $\mathrm{Aut}(X,\mathcal{B},\mu )$.
In other words, $\mathrm{A}(S)$ corresponds to the set of probability measure
preserving (p.m.p.) actions of the free group $\F_{S}$ on 
$(X,\mathcal{B},\mu)$. The metric $\delta $ extends to $\mathrm{A}(S)$ and defines the weak
topology on it. For more details see the book of Kechris \cite{Kechrisbook}.
Note that the space $\mathrm{A}(S)$ contains all p.m.p. actions of groups generated by 
$\left\vert S\right\vert $ elements, where words in $\F_S$ that evaluate to 1 act trivially.  

\vi
Inspired by a classical concept in representation theory, Kechris \cite{Kechrisbook} introduced the weak containment relation on $A(S)$ as follows.
\begin{defin}\label{weakeq}
 Let $f,g\in A(S)$ be two $\F_S$-actions.  Then $f$ {\bf  weakly
contains $g$ }($f\succeq g$) if for any $n\geq 2$, Borel partition $C:X\rightarrow
\{1,2,\dots ,n\}$, a finite set $F\subset \F_{S} $ and $\varepsilon >0$
there exists a Borel partition $D:X\rightarrow \{1,2,\dots ,n\}$ such that 
\begin{equation*}
|\mu (f_{\gamma }D_{i}\cap D_{j})-\mu (g_{\gamma }C_{i}\cap C_{j})|\leq
\varepsilon \quad (1\leq i,j\leq n,\gamma \in F)\,,
\end{equation*}%
where $C_{k}=C^{-1}(k)$ and $D_{k}=D^{-1}(k)$.
\end{defin} 
\noindent This means that the way $g$
acts on finite partitions of the standard Lebesgue space can be simulated by 
$f$ with arbitrarily small error. The actions $f$ and $g$ are {\bf weakly equivalent}, $f\sim g$
if they weakly contain each other. If $\Gamma$ is a countable group, the weak equivalence of probability 
measure preserving actions of $\Gamma$ can be defined analogously to the weak equivalence of $\F_S$-actions.

\vi
One can immediately observe that if the actions $f$ and $g$ are conjugates then they
are weakly equivalent. However, the converse is far from being true, e.g., all
essentially free p.m.p. actions of the integers are pairwise weakly
equivalent. Nevertheless, Kechris (\cite{Kechrisbook}, Proposition 10.1) gave an equivalent definition for weak
containment that involves the notion of conjugacy.
 Let $f,g\in \mathrm{A}(S)$. Then $f$ weakly
contains $g$ ($f\succeq g$) if $C(g)\subseteq \overline{C(f)}$ where the
closure is in the weak topology (\cite{Kechrisbook}, Proposition 10.1). Here, $C(g)$ denotes the conjugates of $g$ in $\mathrm{A}(S)$.
\vi
The aim of this paper is to introduce and study the notion of {\bf partition metric}, a natural
pseudometric on $\mathrm{A}(S)$ such that the zero classes are exactly the
weak equivalence classes and the metric identification of $\mathrm{A}(S)$
with respect to this pseudometric is compact. In the following we give a
condensed description of how to define the partition metric. For details see
Section \ref{sect2}.
\vi
One can extract the local structure of an action $f\in \mathrm{A}(S)$ by taking the stabilizer of a $\mu $-random point in $X$. 
This gives a probability measure on the space of
subgroups of $\F_{S}$ invariant under conjugation by $\F_{S}$. Such
measures are called {\bf invariant random subgroups (IRS)}. The term 
was
coined in \cite{agv}, where it is proved that every IRS arises as the
stabilizer of a p.m.p. action. The IRS of an action encodes information about its freeness but omits many other details.  For instance, all free actions have trivial IRSs. 
\begin{propo} \label{irspropo}
Weakly equivalent actions have the same IRS.
\end{propo}
\noindent
 For an IRS $\lambda $, let $\mathrm{A}(S,\lambda )$ (the 
{\bf fiber} of $\lambda $) be the set of actions in $\mathrm{A}(S)$ with
IRS  $\lambda $.
\vi
The notion of symbolic dynamics can be extended from free actions to an arbitrary element $f\in \mathrm{A}(S)$ as
follows. A Borel partition $C:X\rightarrow \{1,\ldots ,k\}$ defines a
colored
random rooted graph. This graph is constructed as the $k$-vertex colored Schreier
graph of the action on the orbit of a $\mu $-random point of $X$. Applying
Hausdorff distance on the set of $k$-Borel partition of $X$ leads to the
{\bf partition metric} $\mathrm{pd}(f,g)$. 
\vi
This introduces a third equivalent definition of weak equivalence, complementing those given by Kechris. Specifically, two actions \( f \) and \( g \) are weakly equivalent if their partition distance equals zero.

\begin{theorem}
\label{keytheorem} The function $\mathrm{pd}$ is a pseudometric on $\mathrm{A%
}(S)$.  The zero classes of \(\mathrm{pd} \) correspond exactly to the weak equivalence classes. Also, the metric identification of \( A(S) \) with this pseudometric is compact. Moreover, the fiber $A(S,\lambda)$
 of any invariant random
subgroup $\lambda$ is compact as well.
\end{theorem}
\noindent
This metric identification is referred to as the space of actions modulo weak equivalence. The type of convergence within this space is called {\bf local-global
convergence}.
\noindent
The notions of partition metric and local-global convergence come from
graph convergence theory. 
Partition metric has been introduced by Bollob\'{a}s and Riordan in that
setting \cite{Bollobas}. Hatami, Lov\'{a}sz and Szegedy \cite{HLS} later developed the concept of local-global convergence for sequences of finite graphs with bounded degrees, demonstrating that graphings can serve as limit objects. The graph theoretic analogue of an IRS is a unimodular random
rooted graph that has been introduced by Aldous and Lyons in \cite{AL}.
\vi
 An intriguing property of the partition metric is that the weak equivalence class of an action approximately satisfying a combinatorial property  always contains an action satisfying it exactly. We demonstrate this phenomenon on the measurable von 
Neumann problem. 
\noindent
A conjecture attributed to John von Neumann posits that any non-amenable group contains a free subgroup on two generators. However, this conjecture is not universally valid, as non-amenable torsion groups exist.
\cite{Ols}. Nevertheless, the conjecture remains unresolved in the measurable setting, as stated below (\cite{Gab}).

\begin{problem}[Gaboriau-Lyons]
\label{ize}Is it true that for any free p.m.p. action of a countable
non-amenable group $\Gamma $ there exists a free p.m.p. action of $\F_{2}$ on
the same space, such that for $\mu $-almost all $x\in X$, the $x$-orbits
satisfy $x^{\F_{2}}\subseteq x^{\Gamma }$?
\end{problem}
\noindent
The conjecture was settled in the affirmative for large enough Bernoulli
actions by Gaboriau and Lyons \cite{Gab}. 
Using their result, we prove
the following theorem.

\begin{theorem}
\label{negyes} \label{Gab} For any free p.m.p. action of a countable
non-amenable group $\Gamma $ there exists a weakly equivalent action for
which Problem \ref{ize} has an affirmative solution.
\end{theorem}
\noindent
The novel technique introduced in our paper is the use of ultraproducts. This
allows us to extend our results to unitary representations.
\vi
\begin{remark} The first version of this paper was uploaded to arXiv more than ten years ago. Since then, it has been cited and used in several research papers. Partly because of that, we decided not to change the mathematical content in this version. We hope that we made the exposition more clear. 
\end{remark}
\vi The paper is built as follows. In Section \ref{sect2} we define the basic notions and prove some lemmas. In Section \ref{ultratech} we introduce the ultraproduct technique needed. We prove compactness of the space of actions in Section \ref{compactness}. This will lead to the proof of Theorem \ref{negyes} in Section 
\ref{negyes}. In Section \ref{techno} we prove Proposition \ref{irspropo} and establish the equivalence of Kechris' notion of weak equivalence and our version using symbolic dynamics. Finally, in the short Section \ref{vege} we finish the proof of Theorem \ref{keytheorem}. Note that in a follow-up paper Tucker-Drob \cite{TD} defined another pseudo-metric on the space of actions using Kechris' original notion and showed its compactness (using the ultraproduct technique). 

\section{Preliminaries\label{sect2}}
\noindent
In this section we provide the basic definitions and some relevant lemmas.
\vskip 0.2in
\noindent \textbf{Invariant random subgroups and Schreier graphs.} Let $S$ be a finite set and let
the group $\F_S$ act 
transitively by permutations on the pointed countable set $X$. We define the 
\emph{Schreier graph} of this action as follows: the vertex set is $X$ and
for each $s\in S$ and vertex $x$, there is an $s$-labeled edge going from $x$
to $x^{s}$. Let us root the Schreier graph at the distinguished point of $X$. Then
the Schreier graph is a rooted, connected, edge-labeled graph. We identify
Schreier graphs that are isomorphic as rooted, edge labeled graphs. A
particular case is when $H$ is a subgroup of $\F_S $ and the action is the
right coset action: we denote the corresponding Schreier graph by $\mathrm{%
Sch}(\F_S /H,S)$, rooted at $H$. It is easy to see that every Schreier
graph can be obtained this way and that $H$ can be obtained by evaluating 
all the
returning walks in the graph using the edge labels.
\vi
For an abstract alphabet $S$ let $SC(S)$ denote the set of isomorphism
classes of Schreier graphs for the free group $\F_{S}$. An external way to
get an element of $SC(S)$ is to take any $2\left\vert S\right\vert $-regular
graph and then label the directed edges by $S\cup S^{-1}$ such that the
following hold:

\begin{enumerate}
\item for every vertex $x$ and every $s\in S\cup S^{-1}$, there is exactly
one $s$-labeled edge leaving and arriving to $x$;

\item for every directed edge, its label is the formal inverse of the label
of the reverted edge.
\end{enumerate}
\noindent
For two rooted Schreier graphs $G_{1}$ and $G_{2}$, let the distance $%
d(G_{1},G_{2})=\frac{1}{r+1}$, where $r$ is the maximal integer such that the $r$-balls
around the root of $G_{1}$ and $G_{2}$ are isomorphic. The metric $d$ turns $%
SC(S)$ into a totally disconnected, compact space. The group $\F_{S}$ acts on $%
SC(S)$ continuously by moving the root along the path that represents the
acting word.
\vi
Let $SC_{k}(S)$ denote the set of rooted Schreier graphs together with a $k$%
-vertex coloring. We can define the metric similarly as for ordinary
Schreier graphs, just that we consider vertex-colored isomorphisms of rooted
balls in the definition of $r$. Again, this metric turns $SC_k(S)$ to a
totally disconnected, compact space and $\F_{S}$ acts on $SC_k(S)$ continuously
by moving the root. Clearly, the color-forgetting map $SC_{k}(S)\rightarrow
SC(S)$ is a continuous $\F_{S}$-equivariant surjection.
\vi
Let $\mathrm{Sub}(\F_S )$ denote the set of subgroups of $\F_S $. We can
endow $\mathrm{Sub}(\F_S )$ with the topology inherited from the product
topology on the set of subsets of $\F_S $: $\{0,1\}^{\F_S}$. This turns $\mathrm{Sub}(\F_S
)$ into a compact space. The group $\F_S $ acts on $\mathrm{Sub}(\F_S )$
continuously by conjugation. A random subgroup of $\F_S $ is called an 
\emph{invariant random subgroup} (IRS) if its distribution is a Borel
measure invariant under the conjugation action. The name IRS was coined in \cite{agv}.
\vi
For $f\in \mathrm{A}(S)$ let the {\bf type of }$f$ be $\mathrm{Stab}%
_{\F_S }(x)$ where $x$ is a uniform $\mu $-random point in $X$. It is easy
to see that the type is an IRS of $\F_S $. In \cite{agv} it is proved that
every IRS arises as the type of a p.m.p. action. For an IRS $\lambda $ let $%
\mathrm{A}(S,\lambda )$ the {\bf fiber of }$\lambda $ be the set of actions
in $\mathrm{A}(S)$ with type $\lambda $. Another way to look at the type of
an action $f\in \mathrm{A}(S)$ is to consider the Schreier graph of the
action of $\F_{S}$ on the orbit of a uniform $\mu $-random point in $X$,
rooted at $x$. From this point of view, the type is a Borel probability
distribution on $SC(S)$ that is invariant under moving the root. This
identification matches with the canonical bijection between $SC(S)$ and $%
\mathrm{Sub}(\F_{S})$, so there is no ambiguity. \bigskip

\noindent \textbf{The partition metric.} Let $\mathrm{U}(S)$ and$\ \mathrm{U}%
_{k}(S)$ denote the set of $\F_{S}$-invariant Borel probability distributions
on $\mathrm{SC}(S)$ and $\mathrm{SC}_{k}(S)$, respectively. Both sets are compact and metrizable under the weak topology. Let $d_{S}$ and $d_{S,k}$ denote the metrics on $\mathrm{SC}(S)$ and $\mathrm{SC}_{k}(S)$, respectively, which define the weak topology on these sets.  Let us endow the set of compact subsets of $%
\mathrm{U}_{k}(S)$, $\mathrm{Comp}(U_k(S))$ with the Hausdorff metric $d_{Haus,S,k}$. It is important to note that
both $\mathrm U_k(S)$ and $\mbox{Comp}(\mathrm{U}_k(S))$ are compact metric spaces with respect to
the metrics defined above. Note that the Hausdorff-metric $d_{Haus,S,k}$ depends on the choice of the metric $d_{S,k}$. However, the topology of $\mathrm{Comp}(U_k(S))$ does not. 
\vi
For $f\in \mathrm{A}(S)$ and a Borel partition $C:X\rightarrow \{1,\ldots
,k\}$ let $\Psi^C_f:X\to SC_{k}(S)$ be the Borel map, where
$\Psi^C_f(x)$ is the $k$-vertex colored Schreier graph of the action of $%
\F_{S}$ on the orbit of $x$ (with $x$ as the root). Then the push-forward measure $(\Psi^C_f)_* (\mu)$
is an element of $\mathrm{U}_{k}(S)$.
\begin{defin} 
By considering all possible $k$
-Borel partitions $C$ of $X$, the closure of $\cup_C(\Psi^C_f)_* (\mu)$ 
 gives us a compact subset
 $H^k_S(f)\in \mbox{Comp}(U_k(S))$. This set is
called the {\bf global }$k${\bf -type of }$f$.
\end{defin}

\begin{defin}
Let $f,g\in \mathrm{A}(S)$ be two actions. The $k$-partition pseudodistance $%
\mathrm{pd}_{k}(f,g)$ is defined as follows.
$$pd_k(f,g)=d_{Haus,S,k}\left((\Psi^C_f)_* (\mu), (\Psi^C_g)_* (\mu) \right). $$
\end{defin}
\begin{defin}
\noindent The 
\emph{partition pseudometric}, denoted as $\mathrm{pd}(f,g)$  is given by
\begin{equation*}
\mathrm{pd}(f,g)=\sum_{k=1}^{\infty }\frac{1}{2^{k}}\mathrm{pd}_{k}(f,g)\text{.}
\end{equation*}
\end{defin}
\noindent
Thus, the actions $f$ are encoded with the element $\prod^\infty_{k=1} H^k_S(f)$
of the compact ``codespace'' $\prod^\infty_{k=1} \mbox{Comp}( U_k(S))$.
The main theorem entails that if for a sequence of actions
$\{f\}^\infty_{n=1}$, the associated codes are convergent to an element $s$ of the space
$\prod^\infty_{k=1} \mbox{Comp}( U_k(S))$,
then we have an action $f\in A(S)$ such that $s$ is the code of $f$. 
\vi
\noindent \textbf{Unitary representations.} Let $\Gamma $ be a countable
group and $\alpha ,\beta :\Gamma \rightarrow U(H)$ be unitary
representations of $\Gamma $ on a complex separable Hilbert space $H$. We
say that $\beta $ weakly contains (in the sense of Zimmer) \cite{Kechrisbook}
$\alpha $ if for any finite orthonormal system $v_{1},v_{2},\dots ,v_{n}$ in 
$H$, a finite set $F\subset \Gamma $, and a real number $\varepsilon >0$
there exists an orthonormal system $w_{1},w_{2},\dots ,w_{n}$ such that for
any $1\leq i,j,\leq n$ and $\gamma \in F$ 
\begin{equation*}
|\langle \alpha (\gamma )v_{i},v_{j}\rangle -\langle \beta (\gamma
)w_{i},w_{j}\rangle |<\epsilon \,.
\end{equation*}%
We say that two representations are weakly equivalent if they weakly contain
each other in the sense of Zimmer. Note that the original definition of weak
containment and weak containment in the sense of Zimmer are slightly
different (see Appendix H \cite{Kechrisbook}). In our paper, weak containment
always means weak containment in the sense of Zimmer. We say that a unitary
representation $\alpha$ \emph{contains} $\beta$, if $\beta$ is isomorphic
to a subrepresentation of $\alpha$.
Now fix an unitary
representation $\alpha $. Let us consider the countable set of pairs $(F,n)$%
, where $F\subset \Gamma $ is a finite set and $n\geq 1$ is a natural
number. For any such pair we have a product set $D^{n^{2}\times |F|}=C_{F,n}$,
where $D$ is the unit disc of the complex plane.
Let $K_{F,n}(\alpha )$ be the closure of the set 
\begin{equation*} \{\oplus_{\gamma\in F}\oplus_{1\leq i,j\leq n}
\langle \alpha(\gamma) v_i,v_j\rangle\,\mid\,(v_1,v_2,\dots,v_n) \,
\mbox{\,is an orth. syst.}\}
\end{equation*}%
in $C_{F,n}$. Again, we associate a closed subset $Q(\alpha )$ of a compact
product set to a representation by 
\begin{equation*}
Q(\alpha )=\prod_{F,n}K_{F,n}(\alpha )\subset \prod_{F,n}C_{F,n}\,.
\end{equation*}%
The following lemma is straightforward.

\begin{lemma} The unitary representation
$\alpha$ is weakly equivalent to $\beta$ if and only if $Q(\alpha)=Q(\beta)$.
\end{lemma}
\noindent
We shall prove the following analogue of Theorem \ref{keytheorem}.

\begin{theorem}\label{harmas}
\label{masodik} The image of $Q$ is compact.
\end{theorem}

\noindent \textbf{Measure algebras.} In \cite{Kechrisbook}, the author uses
the measure algebra formalism instead of measure spaces. Similarly, we find the use of measure algebras advantageous in our proofs. Hence, in this
subsection we list some well-known facts about measure algebras and group
actions of measure algebras. A measure algebra $\mathcal{M}$ is a Boolean
algebra with a finitely additive measure $\mu $ that is complete metric
space with respect to the distance $d(A,B)=\mu (A\triangle B)$. If $(X,\mu )$
is a Lebesgue probability space, then the equivalence classes of Borel sets
(two sets are equivalent if their symmetric distance has measure zero) form
a measure algebra, the Lebesgue algebra. Any separable atomless measure
algebra is in fact isomorphic to the Lebesgue algebra. In general, if $(X,%
\mathcal{A},\mu )$ is a measure space with a sigma-algebra, then $\mathcal{M}%
(X,\mu )$ denotes the associated measure algebra. Let $\alpha :\mathcal{M}%
(X,\mu )\rightarrow \mathcal{M}(Y,\nu )$ be an injective (measure preserving)
homomorphism between measure algebras. 
Then there exists a surjective Borel map $\Phi _{\alpha
}:Y\rightarrow X$ such that for any Borel set
$A\subseteq X$, $\overline{\Phi _{\alpha
}^{-1}(A)}=\alpha (\overline{A})$, where $\overline{A}$ denotes the element
of the measure algebra represented by the set $A$. 

\begin{defin}\label{asso}
Let $\psi :\mathbb{F}%
_{S}\rightarrow Aut(\mathcal{M}(X,\mu ))$ be a representation of $\mathbb{F}%
_{S}$ by measure algebra automorphisms. Then there exists $f_{\psi }\in
A(S)$ such that for any $\gamma \in \mathbb{F}_S $ and Borel set $A$, $\psi (\gamma )
(\overline{A})=%
\overline{f_{\psi}(\gamma)(A)}\,.$  The action  $f_\psi$ is called the action
{\it associated} with $\psi$.
\end{defin}
\noindent
If $\psi:\mathbb{F}%
_{S}\rightarrow Aut(\mathcal{M}(X,\mu ))$ and $\phi :\mathbb{F}%
_{S}\rightarrow Aut(\mathcal{M}(Y,\nu ))$  are representations and $\alpha :%
\mathcal{M}(X,\mu )\rightarrow \mathcal{M}(Y,\nu )$ is an injective
 measure preserving
isomorphism commuting with the representations, then the associated map $%
\Phi _{\alpha }$ commutes (up to null sets) with the associated actions $f_{\psi },f_{\phi
}\in A(S)$. Then we say that $f_{\phi}$ \emph{contains} $f_{\psi}$.

\section{The ultraproduct technique}\label{ultratech}

\label{sect3}
The main technical tools in our paper are the ultraproducts of actions and representations.
In this section we briefly recall the construction of ultrapowers of
probability measure spaces from \cite{ESZ}. Let $(X,\mu)$ be a standard
Borel probability measure space and $\omega$ be a nonprincipal ultrafilter.
Let $\lim_\omega$ be the associated ultralimit $\lim_\omega:l^\infty\to%
\mathbb{R}$. The ultrapower of the set $X$ is defined the following way.
\noindent Let $\widetilde{X}=\prod^\infty_{i=1} X_i$, where each $X_i$ is a
copy of our $X$, equipped with the atomless probability measure $\mu$, that
we denote by $\mu_i$ to avoid confusion. We say that $\widetilde{p}%
=\{p_i\}^\infty_{i=1}, \widetilde{q}=\{q_i\}^\infty_{i=1}\in \widetilde{X}$
are equivalent, $\widetilde{p}\sim\widetilde{q}$, if 
\begin{equation*}
\{i\in \mathbb{N}\mid p_i=q_i\}\in \omega\,.
\end{equation*}
Define $\mathbf{X}:=\widetilde{X}/\sim$. Now let $\mbox{$\cal P$}(X_i)$
denote the Boolean algebra of Borel subsets of $X_i$, with measure $\mu_i$.
Then let $\widetilde{\mbox{$\cal P$}}%
=\prod^\infty_{i=1}\mbox{$\cal P$}(X_i)$ and $\mbox{$\cal P$}=\widetilde{%
\mbox{$\cal P$}}/I$, where $I$ is the ideal of elements $\{A_i\}^%
\infty_{i=1} $ such that $\{i\in \mathbb{N}\mid A_i=\emptyset\}\in \omega\,.$
Notice that elements of $\mbox{$\cal P$}$ can be identified with certain
subsets of $\mathbf{X}$: If 
\begin{equation*}
\overline{p}=[\{p_i\}^\infty_{i=1}]\in \mathbf{X}\,\,\mbox{and}\,\, 
\overline{A}= [\{A_i\}^\infty_{i=1}]\in \mbox{$\cal P$}
\end{equation*}
then $\overline{p}\in \overline{A}$ if and only if $\{i\in \mathbb{N}\mid p_i\in
A_i\}\in \omega\,.$ Clearly, for $\overline{A}= [\{A_i\}^\infty_{i=1}]$ and $%
\overline{B}= [\{B_i\}^\infty_{i=1}]$ the following hold:

\begin{itemize}
\item $\overline{A}^c=[\{A^c_i\}^\infty_{i=1}]\,,$

\item $\overline{A}\cup \overline{B}=[\{A_i\cup B_i\}^\infty_{i=1}]\,,$

\item $\overline{A}\cap \overline{B}=[\{A_i\cap B_i\}^\infty_{i=1}]\,.$
\end{itemize}
\noindent
Thus, $\mbox{$\cal P$}$ forms a Boolean algebra on $\mathbf{X}$. An important subalgebra of \( P \), denoted \( P' \), is associated with sequences where, for a Borel set \( A \in X \), \( A_i = A \) for all \( i \).
 Clearly, the Boolean algebra $\mbox{$\cal P$}^{\prime }$ is
isomorphic to the Boolean algebra of Borel sets of $X$. Define $\overline{%
\mu}(\overline{A})=\lim_\omega \mu_i(A_i)$. Then $\overline{\mu}:%
\mbox{$\cal
P$}\to\mathbb{R}$ is a finitely additive probability measure. We will call $%
\overline{A}= [\{A_i\}^\infty_{i=1}]$ the \textit{ultraproduct} of the sets $%
\{A_i\}^\infty_{i=1}$.

\begin{defin}
$N\subseteq \mathbf{X}$ is a \textit{nullset} if for any $\varepsilon>0$
there exists a set $\overline{A_\varepsilon}\in\mbox{$\cal P$}$ such that $%
N\subseteq \overline{A_\varepsilon}$ and $\overline{\mu}(\overline{%
A_\varepsilon})\leq \varepsilon$. The set of nullsets is denoted by $%
\mbox{$\cal N$}$.
\end{defin}

\begin{defin}
We call $B\subseteq \mathbf{X}$ a \textit{measurable set} if there exists $%
\widetilde{B}\in \mbox{$\cal P$}$ such that $B\triangle \widetilde{B}\in %
\mbox{$\cal N$}$.
\end{defin}

\begin{propo}
\cite[Proposition 2.2]{ESZ} The measurable sets form a $\sigma$-algebra $%
\mbox{$\cal B$}_\omega$, and $\overline{\mu}(B)= \overline{\mu}(\widetilde{B}%
) $ defines a probability measure on $\mbox{$\cal B$}_\omega$. We denote
this measure space by $( \mathbf{X}, \overline{\mu})$.
\end{propo}
\noindent
It is important to
note that the measure algebra of this space is not separable.
\vskip 0.1in
\noindent
Let $\{ f^i\}^\infty_{i=1}\subset A(S)$. The ultraproduct of these actions $%
\overline{f}$ is defined the following way. 
\begin{equation*}
\overline{f}_\gamma([\{p_i\}^\infty_{i=1}]=
[\{f^i_\gamma(p_i)\}^\infty_{i=1}]\,.
\end{equation*}
This way we defined a measure- preserving action of $\mathbb{F}_S$ on the
ultraproduct space. If \( f^i = f \) for all \( i \), then \( \overline{f}=f_\omega \) is called the ultrapower of \( f \).
 A $\sigma$-algebra
$\mbox{$\cal B$}'_\omega\subset \mbox{$\cal B$}_\omega$ is 
$\mathbb{F}_S$-invariant if, for any $\gamma\in \mathbb{F}_S$ and
$U\in\mbox{$\cal B$}'_\omega, f_\omega(\gamma)(U)\in \mbox{$\cal B$}'_\omega$.

\begin{propo}
\label{p42} Let $f\in A(S)$ and $f_\omega$ be its ultrapower. Let $%
\mbox{$\cal B$}_\omega^{\prime }$ be a $\F_S$-invariant separable subalgebra
of $\mbox{$\cal B$}_\omega$ on ${\bf X} $
containing the algebra $\mbox{$\cal P$}^{\prime}.$ 
Then the associated $\mathbb{F}_S$-action $g\in A(S)$ 
(see Definition \ref{asso})
is weakly equivalent to $f$.
\end{propo}

\smallskip\noindent\textit{Proof.} The measure algebra $\mathcal{M}(\mathbf{X%
},\mbox{$\cal P$}^{\prime })$ is isomorphic to the measure algebra $\mathcal{%
M}(X,\mu)$, hence $g$ contains $f$.

Now let $A:\mathbf{X}\to \{1,2,\dots, k\}$ be a measurable partition of 
$\mathbf{X}$ such that $A^{-1}(i)\in \cal{B}'_\omega$ and let $%
V_1,V_2,\dots, V_k$ be elements of $\cal{P}$ such that
$\overline\mu(A^{-1}(i)\triangle V_i)=0$ for any $1\leq i \leq k$. Then we have a sequence of Borel
partitions $\{X=V^i_1\cup V^i_2\dots\cup V^i_k\}^\infty_{i=1}$ such that $%
[\{V^i_j\}_{i=1}^\infty]=V_j$, for any $1\leq j \leq k\,.$ By definition,
for any $\varepsilon>0$ and $\gamma\in \mathbb{F}_S$ the set 
\begin{equation*}
A^\gamma_{j,l}:=\{i\,\mid\, |\mu(V^i_j\cap f_\gamma (V^i_l))-\overline{\mu}%
(V_j\cap g_\gamma V_l)|<\varepsilon\}
\end{equation*}
is in the ultrafilter $\omega$. Therefore the action $f$ weakly contains the
action $g$.\qed

\begin{propo}
\label{kulcs} Let $f$ and $f_\omega$ be as described earlier. Suppose $h\in A(S)$ is weakly contained in $%
f$. Then, there exists a $\mathbb{F}_S$-invariant
separable subalgebra $\mbox{$\cal B$}_\omega^{\prime }$ of $\mbox{$\cal B$}_\omega$ containing $\mbox{$\cal P$}%
^{\prime }$, such that the associated $%
\mathbb{F}_S$-action $g$ contains $h$.
\end{propo}

\smallskip\noindent\textit{Proof.} Identify $X$ with the product
space $\prod \{0,1\}$, equipped with the usual product measure $\mu$. If $s$ is a $0-1$%
-string of length $k$ let $A_s$ be the elements of $X$ starting with the
string $s$. Let $St(n)$ denote the set of strings of length $n$.
 Let $\gamma_1,\gamma_2,\dots$ be an enumeration of the elements
of $\mathbb{F}_S$. By our assumption, for any $n\geq 1$ there exists a
Borel-partition of $X$ into $2^n$ pieces 
\begin{equation*}
\bigcup_{s_i\in St(n)} B^n_{s_i}=X,
\end{equation*}
such that for any $1\leq i\leq k$ and strings $s_a$, $s_b$ we have that
\begin{equation*}
|\mu(B^n_{s_a}\cap f_{\gamma_i} B^n_{s_b})-\mu(A_{s_a}\cap h_{\gamma_i}
A_{s_b})|< \frac{1}{10^n}\,.
\end{equation*}
Notice that for $k\leq n$, $B^n_s$ is well-defined if $s$ is a string of
length $k$. Simply let 
\begin{equation*}
B^n_s=\bigcup_{s_i\,\mbox{ starts with}\, s} B^n_{s_i}\,.
\end{equation*}
Observe that 
\begin{equation*}
|\mu(B^n_{s}\cap f_{\gamma_i} B^n_{s^{\prime }})-\mu(A_{s}\cap h_{\gamma_i}
A_{s^{\prime }})|< \frac{1}{2^n}\,,
\end{equation*}
if the strings $s$ and $s^{\prime }$ have length not greater than $n$.

Let $\overline{B_s}=[\{ B^n_s\}^\infty_{n=1}]$. Then clearly, 
\begin{equation*}
\overline{\mu}(\overline{B_s} \cap f_\omega(\gamma) \overline{B_{s^{\prime
}}})=\mu(A_s\cap h_{\gamma} A_{s^{\prime }})\,,
\end{equation*}
for all strings $s$,$s^{\prime }$ and $\gamma\in\F_S$. Hence the
subalgebra $\mbox{$\cal C$}_\omega$ generated by the sets $\overline{B_s}$
is $\F_s$ -invariant and $\F_s$-equivariantly isomorphic to the measure
algebra of $(X,\mu)$. Therefore if $\mbox{$\cal B$}_\omega^{\prime }$
contains $\mbox{$\cal C$}_\omega$ then the associated action $f_\omega$
contains $h$%
.\qed

\vskip 0.2in \noindent The following corollary was also proved in \cite%
{Kechriscikk2} (Proposition 4.7)

\begin{theorem} \label{corol31}
If $f\in A(S)$ weakly contains $h\in A(S)$ then there exists $g\in A(S)$
that is weakly equivalent to $f$ that contains $h$. 
\end{theorem}

\subsection{The ultraproduct of unitary representations}

Let $H$ be separable, complex Hilbert space and $\alpha_1, \alpha_2,\dots$
be unitary representations of the countable group $\Gamma$. We define the
ultraproduct of the representations the following way. First we recall the
notion of the ultrapower of $H$. Let
$\prod^\infty_{n=1} H$ be the set of bounded sequences in $H$ and let $V\subset \prod^\infty_{n=1} H$ be the
set of vectors $\{v_n\}^\infty_{n=1}$ such that $\lim_\omega\langle
v_n,v_n\rangle=0\,.$ Clearly, $V$ is a subspace of $\prod^\infty_{n=1} H$
with a well-defined inner product on $\prod^\infty_{n=1}H/V=\prod_\omega H$
by 
\begin{equation*}
\langle [\{v_n\}^\infty_{n=1}], [\{w_n\}^\infty_{n=1}]\rangle =\lim_\omega\langle
v_n,w_n\rangle\,,
\end{equation*}
where $[\{v_n\}^\infty_{n=1}]$ denotes the element in $\prod_\omega H$
representing $\{v_n\}^\infty_{n=1}\in \prod^\infty_{n=1}H$.
It is a standard result that $\prod_\omega H$ is a nonseparable Hilbert
space. The ultraproduct action is defined by 
\begin{equation*}
\alpha_\omega(\gamma)(v)=[\{\alpha_n(\gamma)(v_n)\}^\infty_{n=1}]\,.
\end{equation*}
Clearly, $\alpha_\omega$ is a unitary representation of $\Gamma$. Again, we
consider the special case, when $\alpha_n=\alpha$ for all $n\geq 1$. Let $%
\hat{H}\subset \prod_\omega H$ be the subspace consisting of vectors in the
form $[\{v_i\}^\infty_{i=1}]$, where $v_i=v_j$ for any $i,j\geq 1$. Then we
have the following analogue of Proposition \ref{p42}.

\begin{propo}
Let $\hat{H}\subset K\subset \prod_\omega H$ be a separable $\Gamma$%
-invariant subspace. Then the restriction of $\alpha_\omega$ on $K$ is
weakly equivalent to $\alpha$.
\end{propo}

\smallskip\noindent\textit{Proof.} Clearly, $\alpha_\omega$ weakly contains $%
\alpha$. It is enough to show that $\alpha$ weakly contains $\alpha_\omega$.
Let $\underline{v}_1,\underline{v}_2,\dots,\underline{v}_k\in K$, $F\subset
\Gamma$ be a finite set and $\varepsilon>0$, where $\underline{v}%
_i=[\{v^n_i\}^\infty_{n=1}]$. Let 
\begin{equation*}
S_{\gamma,i,j}=\{n\,\mid\, \langle \alpha(\gamma) v^n_i, v^n_j\rangle-
\langle \alpha_\omega(\gamma) 
\underline{v}_i,\underline{v}_j\rangle|<\epsilon\,\}\,.
\end{equation*}
By the definition of the ultraproduct, $S_{\gamma,i,j}\in\omega$. Hence $%
\cap_{\gamma\in F}\cap_{1\leq i,j \leq k} S_{\gamma,i,j}\in\omega$ as well.
Thus the lemma follows. \qed

\vskip 0.2in \noindent Now we prove the analogue of Proposition \ref{kulcs}.

\begin{propo}
Let $\alpha:\Gamma\to U(H)$ be a representation that weakly contains the representation $\delta:\Gamma\to U(H)$. Then there
exists $\beta:\Gamma\to U(H)$ weakly equivalent to $\alpha$ 
that contains $\delta$. 
\end{propo}

\smallskip\noindent\textit{Proof.} Let $\{v_n\}^\infty_{n=1}$ be an
orthonormal basis for $H$. Enumerate the elements of $\Gamma$, $%
\{\gamma_i\}^\infty_{i=1}$. Since $\alpha$ weakly contains $\delta$, there
exists an orthonormal system $w^n_1,w^n_2,\dots, w^n_n$ such that for any $%
1\leq i,j,k\leq n$ 
\begin{equation*}
|\langle \alpha(\gamma_i) w^n_j, w^n_k\rangle-\langle \delta(\gamma_i)
v_j,v_k\rangle|<\frac{1}{2^k}\,.
\end{equation*}
Let $\underline{w}_j=[\{w^n_j\}^\infty_{n=1}]\in \prod_\omega H\,,$
where $w^n_j=0$, if $n<j$. Then for
any $i,j,k\geq 1$ 
\begin{equation*}
\langle \alpha_\omega(\gamma_i) \underline{w}_j,\underline{w}_k\rangle=
\langle \delta(\gamma_i) v_j,v_k\rangle\,.
\end{equation*}
Hence $\alpha_\omega$ restricted on the $\Gamma$-invariant subspace
generated by $\hat{H}$ and the vectors $\{\underline{w}_j\}^\infty_{j=1}$
contains $\delta$.\qed

\section{Compactness} \label{compactness}
\subsection{The combinatorics of finite balls} \label{combfin}
Define $U^{r,S}$ as the finite family of $r$-balls (up to rooted,
labeled isomorphisms) around the roots of $\mathbb{F}_S$-Schreier graphs, that is, elements of $SC(S)$. We apply the following convention. If $%
x,y\in V(\kappa)$, $\kappa\in U^{r,S}$ and 
\begin{equation*}
d(root(\kappa),x)=d(root(\kappa),y)=r
\end{equation*}
then $x$ and $y$ are not adjacent in $\kappa$. That is, we set vertices on the boundary to be non-adjacent. Let $W^{r,S}$ be the set of
reduced words of length at most $r$ in $\mathbb{F}_S$. For $\kappa\in
U^{r,S} $, we have a partition $P_\kappa$ of $W^{r,S}$ : 
\begin{equation*}
w_1\equiv_{P_\kappa} w_2
\end{equation*}
if $w_1(root(\kappa))=w_2(root(\kappa))\,.$ By our convention, $%
\kappa_1=\kappa_2$ if and only if $P_{\kappa_1}= P_{\kappa_2}\,.$ For $f\in
A(S)$ let $\Psi_f:X\to SC(S)$ be the Borel map, where $\Psi_f(x)$ is the Schreier graph of the action of $\F_S$ on the orbit of $x$ (with $x$ as the root). For a point $x\in X$, we call the $r$-ball around $\Psi_f(x)$ the $%
r $-ball type of $x$ with respect to $f$. If $\kappa\in U^{r,S}$, then $%
T(\kappa)\in SC(S)$ denote the set of Schreier graphs $G$ such that $%
B_r(root(G))\simeq \kappa$. Clearly, $T(\kappa)$ is a clopen set.  Define
\begin{equation*}
T(\kappa,f):= \Psi_f^{-1}(T(\kappa))
\end{equation*}
as the measurable set of points $x\in X$ such that the $r$-ball type of $x$ is $%
\kappa$.
\vi
We also consider labelled balls. Let $U^{r,S,l}$ be the finite set
of all $l$ vertex labelings of the elements of $U^{r,S}$, up to rooted
labeled isomorphisms. Thus we have a map $U^{r,S,l}\to U^{r,S}$ mapping a
vertex labeled graph to the underlying unlabeled graph. Again, for $\tilde{%
\kappa}\in U^{r,S,l}, T(\tilde{\kappa})$ denotes the set of elements $%
\alpha\in SC_l(S)$ such that the $r$-ball around the root of $\alpha$ is
isomorphic to $\tilde{\kappa}$. If $f\in A(S)$ and $D:X\to \{1,2,\dots,l\}$
is a Borel-partition, then the $r$-ball around the root of $\Psi^D_f(x)$ is
called the $(r,l)$-ball type of $x$ with respect to $f$ and $D$. For $\tilde{%
\kappa}\in U^{r,S,l}$, we denote by $T(\tilde{\kappa},f,D)$ the set of vertices $%
x\in X$ with $(r,l)$-type $\tilde{\kappa}$. 

\subsection{The compactness of the global types}
We begin with a simple lemma.

\begin{lemma}
\label{simple} Let $\{f_n\}^\infty_{n=1}\subset A(S)$ be a sequence of
actions, Also, let 
$\{A_n: X\to\{1,2,\dots,l\}\}^\infty_{n=1}$ be a sequence of partitions and $\overline{f}$
respectively $\overline{A}$ be their ultraproducts. Then for any $\tilde{%
\kappa}\in U^{r,S,l}$ 
\begin{equation*}
\lim_\omega \mu(T(\tilde{\kappa},f_n, A_n))=\overline{\mu} (T(\tilde{\kappa},%
\overline{f},\overline{A}))\,.
\end{equation*}
\end{lemma}

\smallskip\noindent\textit{Proof.} Let $x_n\in X_n$ and $\overline{x}%
:=[\{x_n\}^\infty_{n=1}]$. Then, the $(r,l)$-ball type of $x$ with respect to $%
\overline{f}$ is $\tilde{\kappa}$ if and only if the set of numbers $n$ for which
the type of $x_n$ with respect to $f_n$ is $\tilde{\kappa}$, is an element of the
ultrafilter $\omega$. Hence, 
\begin{equation*}
[\{T(\tilde{\kappa}, f_n, A_n)\}^\infty_{n=1}]= T(\tilde{\kappa},\overline{f}%
,\overline{A})
\end{equation*}
and the lemma follows. \qed

\vskip0.2in \noindent Now let $\{f_n\}^\infty_{n=1}$ be as above such that $%
H^k_S(f_n)$ is a Cauchy-sequence in
the space of compact subsets of $U_k(S)$ with the Hausdorff-metric.
\begin{lemma} \label{negyelso}
\begin{equation}  \label{vegsoegyenlet}
\lim_{n\to\infty} H^k_S(f_n)= H^k_S(\overline{f})\,.
\end{equation}
\end{lemma}
\proof
Let $%
p\in H^k_S(\overline{f})$,where $p=(\Psi^{\overline{A}}_{\overline{f}%
})_\star(\mu)\,$ with $\overline{A}=[\{A_n\}^\infty_{n=1}]\,.$ By
Lemma \ref{simple}, it follows that
\begin{equation*}
\lim_\omega(\Psi^{A_n}_{f_n})_\star(\mu)=p\,.
\end{equation*}
Recall that the ultralimit with respect to $\omega$ is well-defined for any
compact metric space. If $\lim_\omega(\Psi^{A_n}_{f_n})_\star(\mu)=p$
then $\lim_{k\to\infty} (\Psi^{A_{n_k}}_{f_{n_k}})_\star(\mu)=p$ for some
subsequence. This proves that for all $k\geq 1$ we have that $\lim_{n\to\infty} H^k_S(f_n)$ contains $H^k_S(%
\overline{f})$. Now let $\{(\Psi^{A_n}_{f_n})_\star(\mu)\}^\infty_{n=1}$ be
a sequence in $U_k(S)$ converging to an element $p\,.$ Then by Lemma \ref{simple} we have that $p\in H_k^S(\overline{f})$.
Therefore $H^k_S(\overline{f})$ contains $\lim_{n\to\infty} H^k_S(f_n)$.
\qed
\begin{propo} \label{compactpropo}
The set $\{\prod^\infty_{k=1} H^k_S(f)\}_{f\in A(S)}$
is a compact subspace of $\prod^\infty_{k=1}\mathrm{Comp}(U_k(S))$.
\end{propo}
\proof
Let $\{f_n\}^\infty_{n=1},\overline{f}$ be as in the previous lemmas. By Lemma \ref{negyelso}, it is enough to prove that there exists $g\in A(S)$ such that for all $k\geq 1$, $H^k_S(g)=H^k_S(\overline{f})$ holds. 
\vi
Let $\mathcal{C}\subset \mbox{$\cal B$}_\omega$ be an $\F_S$%
-invariant separable subalgebra. Then, 
\begin{equation*}
\bigcup_{A^{\prime }}(\Psi^{A^{\prime }}_{\overline{f}})_\star(\mu)\subseteq
\bigcup_{A}(\Psi^{A}_{\overline{f}})_\star(\mu)\,,
\end{equation*}
where the left hand side is taken over the set of all $\mathcal{C}$%
-partitions and the right hand side is taken over the set of all $%
\mbox{$\cal B$}_\omega$-partitions. Therefore, if we pick $\mathcal{C}$ in
such a way that $\bigcup_{A^{\prime }}(\Psi^{A^{\prime }}_{\overline{f}%
})_\star(\mu)$ contains a dense subset of $H^k_S
(\overline{f})$ then $%
H^k_S(g)=H^k_S(\overline{f})$, where $g$ is the action associated to $%
\mathcal{C}$. 
\vskip 0.1in
\noindent {\bf Proof of Theorem \ref{harmas}.} We need to show that the space of representations modulo weak
equivalence is compact. Let $\{\alpha_i\}^\infty_{i=1}$ be a sequence of unitary
representations on the Hilbert space $H$ such that $\{Q(\alpha_i)\}^\infty_{i=1}$ is convergent, that is
for any finite set $F\subset \Gamma$ and $n\geq 1$ $\{K_{F,n}(\alpha_i)\}^%
\infty_{i=1}$ converges to some compact set $L_{F,n}$ in the Hausdorff
metric. We need to prove that there exists a representation $\alpha$ such
that $K_{F,n}(\alpha)=L_{F,n}$, for any $F$ and any $n$. Let $\alpha_\omega$
be the ultraproduct of the $\alpha_i$'s on $\prod_\omega H$. For each $i\geq 1$ and $%
k\geq 1$ pick orthonormal systems $\{v^{i,k,s}_{F,n}\}^n_{s=1}$ such that 
\begin{equation*}
\oplus_{\gamma\in F} \oplus_{1\leq p,q \leq n} \langle \alpha_i(\gamma)
v^{i,k,p}_{F,n}, v^{i,k,q}_{F,n}\rangle=z^k_{F,n}\,
\end{equation*}
\noindent
and the set $\{z^k_{F,n}\}^\infty_{k=1}$ is dense in $L_{F,n}$.
Let $\overline{v}^{k,s}=[\{v^{i,k,s}_{F,n}\}^\infty_{i=1}]$. Let $K$ be the $%
\Gamma$-invariant subspace of $\prod_\omega H$ generated by $\hat{H}$ and
the vectors $\cup_{F,n}\cup^\infty_{k=1}\cup^k_{s=1} \overline{v}%
^{k,s}_{F,n}\,.$ Let $\alpha$ be the restriction of the ultraproduct action 
$\alpha_\omega$ onto 
$K$. By
definition, $L_{F,n}\subseteq K_{F,n}(\alpha)$ for any $F$ and $n$. Now we
prove the converse. Let $x\in K_{F,n}(\alpha)$. Fix a real number $%
\epsilon>0 $. Then there exists $\underline{w}_1,\underline{w}_2,\dots,%
\underline{w}_n\in \prod_\omega H$ such that for any $\gamma\in F$ and $%
1\leq i,j \leq n$ 
\begin{equation*}
|\langle \alpha(\gamma)\underline{w}_i, \underline{w}_j\rangle-x_{%
\gamma,i,j}|<\epsilon\,,
\end{equation*}
where $x_{\gamma,i,j}$ is the coordinate of $x$ associated to the triple $%
(\gamma,i,j)$. By the definition of the ultraproduct, there exist
orthonormal systems $\{t^k_1, t^k_2,\dots, t^k_n\}^\infty_{k=1}\subset H$
such that for any $1\leq p,q\leq n$ 
\begin{equation*}
\lim_\omega\langle \alpha_k(\gamma) t^k_i, t^k_j\rangle=\langle
\alpha(\gamma) \underline{w}_i,\underline{w}_j\rangle\,.
\end{equation*}
Hence we have a subsequence $\{n_k\}^\infty_{k=1}$ such that 
\begin{equation}  \label{egy1}
\lim_{k\to\infty}\langle \alpha_{n_k}(\gamma) t^{n_k}_i, t^{n_k}_j\rangle=
\langle \alpha(\gamma) \underline{w}_i,\underline{w}_j\rangle
\end{equation}
\noindent
for any triple $(\gamma,i,j)$ as above.
Therefore there exists an element $y\in L_{F,n}$ such that each coordinate
of $y$ differs from the corresponding coordinate of $x$ by at most $\epsilon$%
. Consequently, $L_{F,n}=K_{F,n}(\alpha)$\,.\qed

\vskip0.2in \noindent \textbf{Remark:} In \cite[Corollary 4.5]{Kechriscikk2}
the authors prove an interesting compactness result: If $\{a_n\}^%
\infty_{n=1}\subset A(S)$ is a sequence of actions, then there is a
subsequence $n_0< n_1 <n_2\dots$ and $\{b_{n_k}\}^\infty_{k=1}\subset A(S)$
such that $a_{n_k}\sim b_{n_k}$ and $\{b_{n_k}\}^\infty_{k=1}$ converges in $%
A(S)$ in the weak topology.

\vskip 0.1in \noindent The reader may ask what is the relation of this
result to our compactness theorem. In fact the two theorems are independent
as they use a different topology. \cite[Corollary 4.5]{Kechriscikk2} is not
about the compactness of the space of weak equivalence classes since it is
quite possible that the sequence $\{a_n\}^\infty_{n=1}$ converges to an
action $a\in A(S)$ and the sequence $\{b_{n_k}\}^\infty_{k=1}$ converges to
an action $b\in A(S)$ such that $a$ and $b$ are not weakly equivalent.
Indeed, let $a_n=a$ for each $n\geq 1$, where $a$ is a free action of the
free group $\mathbb{F}_S$ that is not weakly equivalent to the Bernoulli
action $b$. Such actions exist e.g. by \cite{Abelek}. By a result of Abert
and Weiss \cite{AW}, $a$ weakly contains $b$. This implies that there exists
a sequence of actions $\{b_n\}^\infty_{n=1}$ such that $b_n$ is equivalent
to $a_n$ and $\{b_n\}^\infty_{n=1}$ converges to $b$ (here we used the 
the definition of weak containment given in the Introduction).

\section{The proof of Theorem \protect\ref{negyes}}

Before starting the proof of the theorem, let us make some remarks.
Weak equivalence of group actions shows some similarity to orbit equivalence
of group actions. It is known that all free ergodic actions of a 
countably infinite
amenable
group are both weakly equivalent and orbit equivalent. By Epstein's theorem 
\cite{Epst} for any non-amenable countable group $\Gamma$ there exist
uncountable many pairwise orbit-inequivalent free actions of $\Gamma$. On
the other hand it is proved in \cite{Abelek} that for several non-amenable
groups there exist uncountable many pairwise weakly-inequivalent actions.
According to Popa's Superrigidity Theorem \cite{Popa} there exist free
actions $\alpha$ of Kazhdan groups $\Gamma$ that are rigid in the
sense, that if an other action is orbit equivalent to $\alpha$ then the two
actions are in fact isomorphic. In \cite{Abelek} it was shown that if two
strongly ergodic profinite actions of a countable group are weakly
equivalent then they are isomorphic. This is however somewhat weaker than
actual rigidity.

\begin{question}
Does there exist a countable group $\Gamma$ with a weakly rigid action ?
\end{question}
\noindent Remark: After this paper was out on arxiv, Robin Tucker-Drob proved that such weakly rigid actions
do not exist \cite{TD}. We leave the question here for historical consistency.
 
\vskip 0.2in \noindent {\bf Proof of Theorem \ref{negyes}.}
Let $\Gamma$ be a countable group and $\alpha: \Gamma\curvearrowright(X,\mu)$
be a free action of $\Gamma$. According to a recent result of Ab\'ert and
Weiss \cite{AW}, $\alpha$ weakly contains all the Bernoulli actions of $\Gamma$.
By  Proposition 4.7 in \cite%
{Kechriscikk2}, there exists an action 
$\beta: \Gamma\curvearrowright(Y,\nu)$ that is weakly
equivalent to $\alpha$ and contains the Bernoulli action $\delta:\Gamma\curvearrowright [0,1]^\Gamma$.
That is there exists a map $\pi$ from $Y$ to $[0,1]^\Gamma $ commuting with
the $\Gamma$-actions. By the theorem of Gaboriau and Lyons, there exists a
free p.m.p. action $\gamma$ of the free group of two generators on $%
[0,1]^\Gamma$ such that for any $t\in \mathbb{F}_2$ and almost all $y\in
[0,1]^\Gamma$, $\gamma(t)(y)=g(y)$ for some $g\in\Gamma$, where $g(y)$ is
the image of $y$ under the Bernoulli action. Now we define the action $%
\gamma^{\prime }$ of $\mathbb{F}_2$ on $(Y,\nu)$ the following way. Let $%
x\in Y$, then $\gamma^{\prime }(t)(x)=\beta(g)(x)$, if $\gamma(t)(\pi(x))=g(%
\pi(x))$. Clearly, $\gamma^{\prime }$ is a free action of $\mathbb{F}_2$
satisfying the condition of Theorem \ref{Gab}.\qed

\section{The weak equivalence notions are equivalent} \label{techno}

The goal of this section is to prove that our definition for weak equivalence
and Kechris's original definition are, in fact, equivalent. That is, the
following theorem holds.
\begin{theorem}\label{tetel5}
The actions $f,g\in A(S)$ are weakly equivalent if and
only if for any $k\geq 1$, $H^k_S(f)=H^k_S(g)$.
\end{theorem}

The ``if'' part is easy. Let $C:X\to\{1,2,\dots,k\}$ be a
Borel-partition. Let $\{C^n:X\to \{1,2,\dots,k\}\}^\infty_{n=1}$ be a
sequence of Borel-partitions such that 
\begin{equation*}
lim_{n\to\infty} ( \Psi^{C^n}_f)_\star (\mu)=(\Psi^C_g)_\star (\mu)
\end{equation*}
in the weak topology. Then for any finite set $F\subset \F_S$ and $%
\varepsilon>0$ 
\begin{equation*}
|\mu(f_\gamma C^n_i\cap C^n_j)-\mu(g_\gamma C_i\cap C_j)|\leq \varepsilon
\quad (1\leq i,j\leq n, \gamma\in F)\,,
\end{equation*}
provided that $n$ is large enough.

\noindent The "only if" part is more complex and will be demonstrated in the next two subsections.
\subsection{The proof of Proposition \ref{irspropo}}
\noindent  
First, we have the following lemma.
\begin{lemma}
Let $f,g\in A(S)$ be $\F_S$-actions. Assume that for every
$r>0$ and $\kappa\in U^{r,S}$, it holds that $\mu(T(\kappa,f))=\mu(T(\kappa,g))$. Then, the IRS's of $f$ and $g$ coincide.
\end{lemma}
\proof
Invariant random subgroups are invariant measures on the compact totally disconnected space 
$\mathrm{Sub}(\F_S)$ which we identified with $SC(S)$ in Section \ref{sect2}. Thus, it is enough to show that for any clopen subset $T(\kappa)\in SC(S)$ 
it holds that $\mu(\Psi_f^{-1}(T(\kappa)))=\mu(\Psi_g^{-1}(T(\kappa)))$ (Section \ref{combfin}). However by definition, $\mu(\Psi_f^{-1}(T(\kappa))=\mu(T(\kappa,f))$, $\mu(\Psi_g^{-1}(T(\kappa))=\mu(T(\kappa,g))$, thus the lemma follows. \qed
\vi
So, we need to prove that if $f,g\in A(S)$ are weakly equivalent then for
every $r>0$ and $\kappa\in U^{r,S}$ we have that $\mu(T(\kappa,f))=\mu(T(\kappa,g))\,.$
The following technical lemma will be crucial. 
\begin{lemma}
\label{key1} Let $f\in A(S)$. Then for any $r>0$ and $\delta>0$ there exists a positive integer $n_{r,\delta,f}$
and
a finite partition $X=\cup_{i=1}^{n_{r,\delta,f}} L_i\cup E_{r,\delta,f}$ with
the following properties.

\begin{enumerate}
\item $\mu(E_{r,\delta,f})<\delta\,.$

\item Each $L_i$ is a subset of $T(\kappa,f)$ for some $\kappa\in U^{r,S}$.
We denote this element $\kappa$ by $\tau(L_i)$.

\item If $w_1,w_2\in W^{r,S}$ and $w_1\not\equiv_{P_{\tau(L_i)}} w_2$ then 
\begin{equation*}
f_{w_1}(L_i)\cap f_{w_2}(L_i)=\emptyset\,.
\end{equation*}
\vi
The equivalence relation $\equiv_{P_\kappa}$ has been defined in Section \ref{combfin}.
Note that by our second condition, if $w_1\equiv_{P_{\tau(L_i)}} w_2$ then 
$f_{w_1}(L_i)=f_{w_2} (L_i)\,.$
\end{enumerate}
\end{lemma}

\smallskip\noindent\textit{Proof.} According to Luzin's theorem, there exists a compact
set $C_\delta\subset X$ such that $\mu(X\backslash C_\delta)<\delta/2\,,$ and
all the coordinates of $f\in A(S)=Aut(X,\mu)^S$ are continuous on $C_\delta$%
. Let $x\in C_\delta$. Define $\lambda(x)$ by 
\begin{equation*}
\lambda(x):=\inf_{w_1,w_2\in W^{r,S}\,,f_{w_1}(x)\neq f_{w_2}(x)}
d_X(f_{w_1}(x), f_{w_2}(x))\,,
\end{equation*}
where $d_X$ is the standard metric on the unit interval $X$. Note that $%
\lambda(x)=0$ if and only if $x$ is a fixed point of the action $f$. Let $%
\chi>0$ be a real number such that 
\begin{equation*}
\mu(x\,\mid\, 0<\lambda(x)<\chi)<\delta/2\,.
\end{equation*}
By uniform continuity, there exists an $\varepsilon>0$ such that if $x,y\in
C_\delta$ and $d_X(x,y)<\varepsilon$ then $d_X(f_w(x),f_w(y))<\chi$ for any $%
w\in W^{r,S}$. Now let $E_{r,\delta,f}:=X\backslash C_\delta \cup
\{x\,\mid\, 0<\lambda(x)<\chi\}\,.$ For $\kappa\in U^{r,S}$, choose an
arbitrary finite partition of $T(\kappa,f)\backslash E_{r,\delta,f}$ by
subsets of diameter less than $\epsilon$. Let $L$ be such a subset, $x,y\in L$
and $w_1\not\equiv_{P_\kappa} w_2\,.$ Then $d_X(f_{w_1}(x),f_{w_2}(x))\geq
\chi$ and $d_X(f_{w_2}(x),f_{w_2}(y))<
\chi$.  Hence, $f_{w_1}(x)\neq f_{w_2}(y)$, that is,  $f_{w_1}(L)$ and $f_{w_2}(L)$ are disjoint subsets. \qed
\vi
We fix the integer $r>0$ and $\delta>0$ until the end of the proof of Theorem \ref{tetel5}.
\vskip0.2in \noindent Now, let us introduce the notion of height for $r$-types.
The set $U^{r,S}$ is an ordered set, where $\kappa\leq \lambda$ if $P_\kappa$ is a
refinement of $P_\lambda$. The height function $h_r:U^{r,S}\to\mathbb{N}$ is
defined the following way. If $\kappa$ is a minimal element, then let $%
h_r(\kappa)=1, \Sigma_r(1)=h^{-1}_r(1)\,.$ If $\kappa$ is a minimal element
in $U^{r,S}\backslash \Sigma_r(1)$, then let $h_r(\kappa)=2,
\Sigma_r(2)=h^{-1}_r(2)\,.$ If $\Sigma_r(1), \Sigma_r(2),\dots,\Sigma_r(k)$
are already defined then let $h_r(\kappa)=k+1$ if $\kappa$ is minimal in the
set $U^{r,S}\backslash \cup^k_{i=1} \Sigma_r(i)$ and let $%
\Sigma_r(k+1)=h_r^{-1}(k+1)\,.$ 
\vi Let $P:=\cup^{n_{r,\delta,f}}_{i=1}
L_i\cup E_{r,\delta,f}$ be a partition of $X$ satisfying the conditions of
Lemma \ref{key1}. Let $\rho>0$. We say that a partition $P^\rho:=\cup^{n_{r,%
\delta,f}}_{i=1} L^\rho_i\cup E^\rho_{r,\delta,f}$ is a $\rho$-simulation of 
$P$ if

\begin{enumerate}
\item $|\mu(L_i)-\mu(L^\rho_i)|<\rho$ for any $1\leq i \leq
n_{r,\delta,f}\,. $

\item $|\mu (f_{w_1}(L_i)\cap f_{w_2}( L_j))- \mu(g_{w_1} (L^\rho_i)\cap
g_{w_2}( L^\rho_j))|<\rho$ for any $1\leq i,j \leq n_{r,\delta,f}$ and $%
w_1,w_2\in W^{r,S}\,.$
\end{enumerate}
\noindent
Note that by weak equivalence, such $\rho$-simulations must exist.

\begin{lemma}
\label{prop3} Let $\kappa\in U^{r,S}$. Then 
\begin{equation}  \label{prop3e1}
\limsup_{\rho\to 0} \mu \left( \cup_{L_i\subset T(\kappa,f)}
L^\rho_i\backslash T(\kappa,g)\right) \leq 3\delta
\end{equation}
\end{lemma}

\smallskip\noindent\textit{Proof.} For $1\leq i \leq n_{r,\delta,f}$ let 
\begin{equation} \label{fontos}
\hat{L}^\rho_i=\{x\in L^\rho_i\,\mid \, g_{w_1}(x) \neq g_{w_2}(x)\,\, %
\mbox{\,if}\,\, w_1\not\equiv_{P_{\tau(L_i)}} w_2\}\,.
\end{equation}
\noindent
Observe that if $x\in \hat{L}^\rho_i$, then the $r$%
-ball type of $x$ with respect to $g$ is less or equal than  the $r$%
-ball type of $x$ with respect to $f$. Also, by
the definition of a $\rho$-simulation 
\begin{equation}  \label{prop3e2}
\lim_{\rho\to 0} \mu(\hat{L}_i^\rho)=\mu(L_i)\,.
\end{equation}

\noindent
Now we need another lemmma to proceed. 

\begin{lemma}
For any $n\geq 1$, 
\begin{equation}  \label{em3}
\left|\mu \left(\bigcup_{\lambda, h_r(\lambda)\leq n} T(\lambda,f)\right)-\mu \left(\bigcup_{\lambda,
h_r(\lambda)\leq n} T(\lambda,g)\right)\right|\leq \delta\,.
\end{equation}
\end{lemma}

\smallskip\noindent\textit{Proof.} By definition, 
\begin{equation*}
\mu \left(\bigcup_{\lambda, h_r(\lambda)\leq n} T(\lambda,f)\right)\leq \sum_{i,
h_r(L_i)\leq n} \mu(L_i) +\delta\,,
\end{equation*}
where $h_r(L_i)$ is defined as $h_r(\lambda)$, if $L_i\subset T(\lambda,f)$.
Also,

\begin{equation*}
\mu\left(\bigcup_{\lambda, h_r(\lambda)\leq n} T(\lambda,g)\right)\geq \sum_{i,
h_r(L_i)\leq n} \mu(\hat{L}^\rho_i)\,,
\end{equation*}
\noindent
by the observation after \eqref{fontos}.
Hence by (\ref{prop3e2}), 
\begin{equation*}
\mu\left(\bigcup_{\lambda, h_r(\lambda)\leq n} T(\lambda,g)\right)\geq \mu\left(\bigcup_{\lambda,
h_r(\lambda)\leq n} T(\lambda,f)\right)-\delta\,.
\end{equation*}
Since $f$ weakly contains $g$ the reverse inequality must hold : 
\begin{equation*}
\mu\left(\bigcup_{\lambda, h_r(\lambda)\leq n} T(\lambda,f)\right)\geq 
\mu\left(\bigcup_{\lambda,
h_r(\lambda)\leq n} T(\lambda,g)\right)-\delta\,.
\end{equation*}
That is 
\begin{equation*}
\left|\mu\left(\bigcup_{\lambda, h_r(\lambda)\leq n} T(\lambda,f)\right)-\mu \left(\bigcup_{\lambda,
h_r(\lambda)\leq n} T(\lambda,g)\right)\right|\leq \delta.\,\quad\qed
\end{equation*}

Now we finish the proof of Lemma \ref{prop3}.
By the definition of the subset $L_i$'s, it follows that
\begin{equation}  \label{em1}
\left|\mu\left(\bigcup_{\lambda, h_r(\lambda) <h_r(\kappa)} T(\lambda,f)\right)-\sum_{L_i, h_r(L_i)<h_r(\kappa)}
\mu(L_i)\right|<\delta
\end{equation}
By the observation after \eqref{fontos},  if $x$ is an element of $\cup_{i, L_i\subset T(\kappa,f)}
\hat{L}^\rho_i \backslash T(\kappa,g)$, then the $r$-ball type of $x$ with respect to $%
g $ is strictly smaller than $h_r(\kappa)$. Also if $x\in \hat{L}^\rho_j$ and $h(L_j)<h_r(\kappa)$
then the $r$-ball type of $x$ with respect to $g$ is strictly smaller than $h_r(\kappa)$,
as well. Therefore, we have that
$$\bigcup_{\lambda, h_r(\lambda)<h_r(\kappa)} T(\lambda,g)\supseteq \bigcup_{j, h_r(L_j)<h_r(\kappa)} \hat{L}^\rho_j\cup \bigcup_{L_i, L_i\subset T(\kappa,f)} \hat{L}^\rho_i \backslash
T(\kappa,g)\,.$$ \noindent
  Thus by (\ref{prop3e2}), if $\rho$ is small enough then 
\begin{equation}  \label{em2}
\mu\left(\bigcup_{\lambda, h_r(\lambda)<h_r(\kappa)} T(\lambda,g)\right)\geq \left(\sum_{j, h_r(L_j)<h_r(\kappa)}
\mu(L_j)\right)-\delta+ \mu\left(\bigcup_{L_i, L_i\subset T(\kappa,f)} L^\rho_i \backslash
T(\kappa,g)\right)\,.
\end{equation}
Adding up the inequalities (\ref{em3}) (for $n=h_r(\kappa)+1$),(\ref{em1}) and (\ref{em2}) we get
the statement of Lemma \ref{prop3}. \qed
\vi
Now we finish the proof of 
Proposition \ref{irspropo}. Since Lemma \ref{prop3} holds for arbitrarily small positive $\delta$ and
$\lim_{\rho\to 0} \mu (L^\rho_i)=\mu(L_i)$ holds for all $1\leq i \leq n_{r,\delta,f}$, we get that $\mu(T(\kappa,g))\leq \mu(T(\kappa,f))\,.$  As the reverse inequality must also hold, we have
the equation
$$\mu(T(\kappa,g))= \mu(T(\kappa,f))\,.$$
\noindent
finishing the proof of  Proposition \ref{irspropo}. \qed
\vskip0.1in \noindent 
\subsection{The end of the proof of Theorem \ref{tetel5}} \noindent The following lemma is a straightforward consequence of Proposition 10.1 from \cite{Kechrisbook}, which was previously mentioned in the Introduction.
We leave the details for the reader to verify. 

\begin{lemma}
\label{exercise} The following statements are equivalent for $f,g\in A(S)$:

\begin{itemize}
\item $f$ weakly contains $g$.

\item For any $m,n,k,l\geq 1$, $\delta>0$, finite set $F\subset\Gamma$, and
partitions $A:X\to\{1,2,\dots,k\}$, $B:X\to\{1,2,\dots,l\}$ there exist
partitions $C:X\to\{1,2,\dots,k\}$, $D:X\to\{1,2,\dots,l\}$ such that for
all $1\leq r_1, r_2,\dots, r_m\leq k$, $1\leq q_1, q_2,\dots, q_n\leq l$, $%
\gamma_1, \gamma_2,\dots, \gamma_m\in F$ and $\delta_1, \delta_2,\dots,
\delta_n\in F$ 
\begin{equation*}
\left|\mu\left(\bigcap^m_{i=1} f_{\gamma_i}(C_{r_i})\cap \bigcap^n_{j=1}
f_{\delta_j}(D_{q_j})\right)- \mu\left(\bigcap^m_{i=1} g_{\gamma_i}(A_{r_i})\cap
\bigcap^n_{j=1} g_{\delta_j}(B_{q_j})\right)\right|<\delta\,.
\end{equation*}
\end{itemize}
\end{lemma}
 \vi Now let $C:X\to \{1,2,\dots, l\}$
be a Borel-partition of $X$. 
It is enough to prove that for any 
$\varepsilon>0$ and $r>0$ there exists a partition $C^{\prime }:X\to
\{1,2,\dots, l\}$ of $X$ such that 
\begin{equation}  \label{kulcsegyenlet}
| (\Psi^{C}_f)_\star (\mu) ( T(\tilde{\kappa}))- (\Psi^{C^{\prime
}}_g)_\star (\mu) ( T(\tilde{\kappa}))|<\varepsilon\,
\end{equation}
holds for all $\tilde{\kappa}\in U^{r,S,l}$. Indeed, it means that $%
H^l_S(f)\subseteq H^l_S(g)$.

\noindent For the rest of the proof we fix an integer $r>0$ and a real number $\delta>0$.
Let $T_\delta=\cup_{i=1}^{n_{r,\delta,f}} L_i\cup E_{r,\delta,f}$ be a
Borel-partition of $X$ as in Lemma \ref{key1}. We say that a pair of
partitions of $X$, $(C^\rho, T^\rho_\delta)$ is a $\rho$-simulation of the
pair $(C,T_\delta)$ if

\begin{itemize}
\item $T^\rho_\delta$ is a $\rho$-simulation of $T_\delta$.

\item For any $\tilde{\kappa}\in U^{r,S,l}$ and $L_i\subset T([\tilde{\kappa}%
],f)$ we have that

\begin{equation}\label{hosszu}
\left|\mu\left(L^\rho_i\cap \bigcap_{w_j\in W^{r,S}} g_{{w_j}^{-1}} (C^\rho_{%
\tilde{\kappa}(w_j)})\right)- \mu\left(L_i\cap \bigcap_{w_j\in W^{r,S}} f_{{%
w_j}^{-1}} (C_{\tilde{\kappa}(w_j)})\right)\right|<\rho\,,
\end{equation}
\end{itemize}
\noindent
where $\tilde{\kappa}(w_j)$ denotes the label of $w_j(root(\tilde{\kappa}))$
in $\tilde{\kappa}$ and $[\tilde{\kappa}]$ denotes the underlying $r$-ball type
of $\tilde{\kappa}$. By Lemma \ref{exercise}, such $\rho$-simulation exists.

\noindent Now by definition,
\begin{equation} \label{ma1}
T([\tilde{\kappa}],f)\cap \bigcap_{w_j\in W^{r,S}} f_{{w_j}^{-1}} (C_{\tilde{\kappa}%
(w_j)})=T(\tilde{\kappa},f,C)\end{equation}
\noindent
Similarly,
\begin{equation} \label{ma2}
T([\tilde{\kappa}],g)\cap \bigcap_{w_j\in W^{r,S}} g_{{w_j}^{-1}} (C_{\tilde{\kappa}%
(w_j)})=T(\tilde{\kappa},g,C^\rho)\end{equation}
\noindent
Hence by the definition of the $L_i$'s, we have that

\begin{equation}\label{ma3}
\left|\mu\left(\bigcup_{L_i, L_i\subset T([\tilde{\kappa}],f)} \left(L_i\cap
\bigcap_{w_j\in W^{r,S}} f_{{w_j}^{-1}} (C_{\tilde{\kappa}%
(w_j)})\right)\right)-(\Psi^C_f)_\star(\mu)(T(\tilde{\kappa}))\right| <\delta
\end{equation}
\vi
Also by \eqref{ma2}, we have that
\begin{equation}\label{ma4}
\begin{split}
\left|\mu\left(\bigcup_{L_i, L_i\subset T([\tilde{\kappa}],f)} \left(L^\rho_i
\bigcap_{w_j\in W^{r,S}} g_{{w_j}^{-1}} (C_{\tilde{\kappa}%
(w_j)})\right)\right)-(\Psi^C_g)_\star(\mu)(T(\tilde{\kappa}))\right| < \\
< \mu
\left(\left(T([\tilde{\kappa}],g)\backslash
\bigcup_{L_i, L_i\subset T([\tilde{\kappa}],f)} L_i^\rho\right)\cup \left(\bigcup_{L_i, L_i\subset T([\tilde{\kappa}],f)} L_i^\rho
\backslash T([\tilde{\kappa}],g)\right)\right)
\end{split}
\end{equation}
Therefore
by Lemma \ref{prop3} and the fact that
$\mu(T([\tilde{\kappa}],f)=
\mu(T([\tilde{\kappa}],g)$, if both $\delta$ and $\rho$ are sufficiently
small then 
$$| (\Psi^{C}_f)_\star (\mu) ( T(\tilde{\kappa}))- (\Psi^{C^{\rho
}}_g)_\star (\mu) ( T(\tilde{\kappa}))|<\varepsilon\,$$
 holds. This finishes the proof of Theorem
\ref{tetel5}. \qed
\section{The proof of Theorem \ref{keytheorem}} \label{vege}
\noindent
Finally, we are able to prove  Theorem \ref{keytheorem}.
By Theorem \ref{tetel5}, the weak equivalence classes correspond exactly to the zero classes of $pd$. Also, by Proposition \ref{compactpropo}, the metric space of the zero classes is compact. The compactness of the fiber $A(S,\lambda)$ follows directly from the definition of the partition pseudometric $pd$.  This finishes the proof of Theorem \ref{keytheorem}. \qed

\end{document}